\documentclass{amsart}
\usepackage{amssymb}
\usepackage{amsfonts}
\usepackage{amsmath}
\usepackage[legalpaper,bookmarks=true,colorlinks=true,linkcolor=blue,citecolor=blue]{hyperref}
\usepackage{graphicx}%
\usepackage{fancyhdr}
\usepackage{color}
\usepackage[mathlines]{lineno}
\usepackage{lscape}
\usepackage{epsfig}
\usepackage{natbib}
\usepackage{geometry}
\usepackage{tgbonum}
\usepackage[]{listings}

\fontfamily{qcr}\selectfont

\newtheorem{theorem}{Theorem}
\theoremstyle{plain}

\newtheorem{proposition}{Proposition}

\numberwithin{equation}{section}

\newcommand{\Bin}{\bigskip \noindent}

\newcommand{\Ni}{\noindent}

\begin{document}
\Large
\title[Extensions of two classical Poisson limit laws to non-stationary independent data]{Extensions of two classical Poisson limit laws to non-stationary independent sequences}

\author{Aladji Babacar Niang} \author{Harouna Sangar\'e} \author{Tchilabalo Abozou Kpanzou}  \author{Gane Samb Lo}  \author{Nafy Ngom}

\begin{abstract} In earlier stages in the introduction to asymptotic methods in probability theory, the weak convergence of sequences $(X_n)_{n\geq 1}$ of Binomial of random variables (\textit{rv}'s) to a Poisson law is classical and easy-to prove. A version of such a result concerning sequences $(Y_n)_{n\geq 1}$ of negative binomial \textit{rv}'s also exists. In both cases, $X_n$ and $Y_n-n$ are by-row sums $S_n[X]$ and $S_n[Y]$ of arrays of Bernoulli \textit{rv}'s and corrected geometric \textit{rv}'s respectively. When considered in the general frame of asymptotic theorems of by-row sums of \textit{rv}'s of arrays, these two simple results in the independent and identically distributed scheme can be generalized to non-stationary data and beyond to non-stationary and dependent data. Further generalizations give interesting results that would not be found by direct methods. In this paper, we focus on generalizations to the non-stationary independent data. Extensions to dependent data will addressed later.\\

\noindent $^{\dag}$ Aladji Babacar Niang\\
LERSTAD, Gaston Berger University, Saint-Louis, S\'en\'egal.\\
Email: niang.aladji-babacar@ugb.edu.sn, aladjibacar93@gmail.com\\
Imhotep Mathematical Center (IMC), imhotepsciences.org\\

\noindent $^{\dag\dag}$ Dr Harouna Sangar\'e\\
Main Affiliation: DER MI, FST, Universit\'e des Sciences, des Techniques et des Technologies de Bamako (USTT-B), Mali.\\
Affiliation : LERSTAD, Universit\'e Gaston Berger (UGB), Saint-Louis, S\'en\'egal.\\
Email : harounasangare@fst-usttb-edu.ml, harouna.sangare@mesrs.ml\\
sangare.harouna@ugb.edu.sn, harounasangareusttb@gmail.com\\

\noindent $^{\dag\dag\dag}$ Dr Tchilabalo Abozou Kpanzou
University of Kara, Kara, Togo\\
Affiliated to LERSTAD, Gaston Berger University, Saint-Louis, Senegal\\
Emails: t.kpanzou@univkara.net; kpanzout@gmail.com\\

\noindent $^{\dag\dag\dag\dag}$ Gane Samb Lo.\\
LERSTAD, Gaston Berger University, Saint-Louis, S\'en\'egal (main affiliation).\newline
LSTA, Pierre and Marie Curie University, Paris VI, France.\newline
AUST - African University of Sciences and Technology, Abuja, Nigeria\\
Imhotep Mathematical Center (IMC), imhotepsciences.org\\
gane-samb.lo@edu.ugb.sn, gslo@aust.edu.ng, ganesamblo@ganesamblo.net\\
Permanent address : 1178 Evanston Dr NW T3P 0J9,Calgary, Alberta, Canada.\\

\noindent $^{\dag\dag\dag\dag\dag}$ Nafy Ngom\\
LERSTAD, Gaston Berger University, Saint-Louis, S\'en\'egal (main affiliation).\newline
fany.ngom@ugb.edu.sn\\

\end{abstract}
\maketitle

\section{Introduction}

\subsection{preliminaries} $ $ \\

\Ni The approximation of a sequence of binomial probability laws $\left(\mathcal{B}(n,p_{n})\right)_{n\geq 1}$ associated to a sequence of \textit{r.v}'s $(Z_n)_{n\geq 1}$ [such that the sequence of probabilities $(p_n)_{n\geq 1}$ converges to zero and $np_n \rightarrow \lambda >0$ as $n\rightarrow +\infty$] to a Poisson law $\mathcal{P}(\lambda)$ is a classical and easy-to-prove result in probability theory. This approximation has very important applications in real-life problems, especially in lack of powerful computers. In this simple case, each $Z_n$ is a sum of $n$ independent and identically distributed (\textit{iid}) Bernoulli $\mathcal{B}(p_n)$-random variables $\left\{(X_{j,n})_{1\leq j \leq n}, \ n\geq 1\right\}$, \textit{i.e.}, $Z_n=X_{1,n}+X_{2,n}+\cdots+X_{n,n}$. When we depart from the identical distributivity assumption, the problem may get more and rapidly complex, even if the independence assumption is still required. The situation becomes more interesting if the random variables $X_{j,n}$ are non-stationary and independent.  \\

\Ni In this note, we aim at giving non trivial generalizations of such results in the frame of the central limit theorem for independent random variables.\\

\Ni Also, there is a \textit{negative}  version of the described result. Indeed, if we call a binomial law as a positive binomial law $\mathcal{PB}(n,p)$, $n\geq 1$, $0<p<1$ in opposition to  a negative binomial law $\mathcal{NB}(n,p)$, we have the two following results concerning positive binomial and negative binomial laws respectively.\\

\Ni Let us make this precision for once: throughout this paper, all limits are meant as $n\rightarrow +\infty$ unless the contrary is specified.\\

\begin{proposition} \label{piidBinomial} Let $(X_{n})_{n\geq 1}$ be a sequence of random variables in some probability space $(\Omega,\mathcal{A}, \mathbb{P})$ such that: \\

\Ni 1) $\forall n\geq 1 $, $X_n \sim \mathcal{B}(n,p_{n})$,  $n\geq 1 $;\\

\Ni 2) $p_{n} \rightarrow 0$ and $np_{n} \rightarrow \lambda \in \mathbb{R}_{+}\setminus\{0\}$ as $n\rightarrow + \infty$.\\

\Ni Then 

$$
X_{n} \rightsquigarrow  \mathcal{P}(\lambda).
$$
\end{proposition}

\Bin Next, we have: \\

\begin{proposition}\label{niidBinomial}
Let $(X_{n})_{n\geq 1}$ be a sequence of random variables in some probability space $(\Omega,\mathcal{A}, \mathbb{P})$ such that: \\

\Ni 1) $\forall n\geq 1 $, $X_n \sim \mathcal{N}\mathcal{B}(n,p_n)$,  $n\geq 1 $; \\

\Ni 2) $(1-p_{n}) \rightarrow 0$ and $n(1-p_{n}) \rightarrow \lambda \in \mathbb{R}_{+}\setminus\{0\}$ as $n\rightarrow + \infty$.\\

\Ni Then 

$$
X_{n}-n \rightsquigarrow \mathcal{P}(\lambda).
$$
\end{proposition}

\Bin \textbf{Remark}. The second result is proved in \cite{ips-wcia-ang}. Although the proof is direct, we do not encounter it in some classical books as \cite{feller1, feller2}, \cite{gutt}, \cite{loeve}. For that reason, we give it in the appendix in page \pageref{niidBinomial-proof}.\\

\Ni Our aim here is to provide non-trivial generalizations of such simple results to non-stationary and independent data. The used methods will allow later further generalizations even with dependent data.\\

\Bin Let us prepare generalizations by transforming both results as sums of random variables. \\

\subsection{The Central limit theorem frame} $ $\\

\Ni It is known that a binomial random variable $X_{n} \sim \mathcal{B}(n,p_{n})$ has the same law than a sum of $n$ \textit{iid} Bernoulli distributed random variables:

$$
X_{n}=^d X_{1,n}+\cdots+X_{n,n},
$$

\Bin where $X_{1,n},\cdots,X_{n,n}$ are independent and follow all the $\mathcal{B}(p_n)$-law. Also $X_n \sim \mathcal{N}\mathcal{B}(n,p_{n})$ has the same law of a sum of $n$ \textit{iid} \textit{r.v}'s: 

$$
X_{n}=^d X^{\ast}_{1,n}+\cdots+X^{\ast}_{n,n},
$$

\Bin where $X^{\ast}_{1,n},\cdots,X^{\ast}_{n,n}$ are independent and each $X^{\ast}_{j,n}$ follows the geometric law $\mathcal{G}(p_{n})$. In the second, we rather use

$$
X_{n}-n=\sum_{i=1}^{n}\left(X^{\ast}_{i,n}-1\right)=:\sum_{i=1}^{n}X_{i,n},
$$

\Bin where the $X_{i,n}$'s are independent and each $X_{i,n}$ follows the law $\mathcal{G}^{\ast}(p_{n})=\mathcal{G}(p_{n})-1$, and such a law is called a corrected geometric law, for convenience. In both cases, we have to study an array

$$
X\equiv \biggr\{ \{X_{k,n}, \ 1 \leq k \leq k(n)\}, \ n\geq 1\biggr\},
$$ 

\Bin of random variables defined in the same probability space $(\Omega, \mathcal{A},\mathbb{P})$ such that here: \\ 

\Ni 1) $\forall n\geq 1 $, $k(n)=n$;\\

\Ni 2) $\forall n\geq 1$, the variables $X_{1,n},\cdots,X_{k(n),n}$ are independent;\\

\Ni 3) The sequence  ${X_{1,n},\cdots,X_{k(n),n}}$ is stationary for $n\geq 1$;\\

\Ni 4) $\forall k \in [1,k(n)]$, $X_{k,n} \sim \mathcal{B}(p_{n})$ or $X_{k,n} \sim \mathcal{G}(p_{n})-1$. \\

\Ni We see that we are in the \textit{CLT} frame and each of points (2) and (3) can be changed to lead to generalizations. Since, we want to generalize the limit binomial laws and negative binomial laws, we keep the same hypotheses on the marginal laws

\begin{equation*}
\forall k \in [1,k(n)], \ X_{k,n} \sim \mathcal{B}(p_{k,n}) \ \ or \ \ X_{k,n} \sim \mathcal{G}(p_{k,n})-1
\end{equation*}

\Bin and try to answer to the questions (Q1) and (Q2) below:

\Bin (Q1) Given an array $X$ of random variables with $k(n)\rightarrow +\infty$ such that the elements of each row are independent and $\mathcal{B}(p_{k,n})$-\textit{r.v}'s, do we still have 

\begin{equation}
S_{n}[X]=\sum_{k=1}^{k(n)}X_{k,n}\rightsquigarrow \mathcal{P}(\lambda), \label{final-result}
\end{equation}

\Bin when some of the assumptions (1), (2) and (3) [but mainly [2] and [3])] are violated, and under what sufficient conditions this should hold?\\

\Ni (Q2) Given an array $X$ of random variables with $k(n)\rightarrow +\infty$ such that the elements of each row are independent and  $\mathcal{G}^{\ast}(p_{k,n})$-\textit{r.v}'s, do we still have \eqref{final-result} when some of the assumptions (1), (2) and (3) [but mainly [2] and [3])] are violated, and under what sufficient conditions this should hold?\\

\Ni Although direct handlings of these questions might be possible, we think that a general and extensible solution resides  in the \textit{CLT} frame, since it will prepare further generalizations for dependent data.\\

\Ni Therefore, we organize the paper as follows. In section \ref{sec-02}, we recall the frame of the \textit{CLT} problem as stated in \cite{loeve}.  In section \ref{sec-03}, we state and prove the results. We conclude the paper by conclusive remarks in Section \ref{sec-04}.  \\

\section{Notation and $G$-\textit{CLT} for summands of independent random variables}\label{sec-02}

\Ni Let us consider the array
  
$$
X\equiv \biggr\{ \{X_{k,n}, \ 1 \leq k \leq k_n=k(n)\}, \ n\geq 1\biggr\},
$$

\Bin of square integrable random variables defined on the same probability space $(\Omega, \mathcal{A}, \mathbb{P})$. We denote $F_{k,n}$ as the cumulative distribution function (\textit{cdf}) of $X_{k,n}$. We also denote by $a_{k,n} = \mathbb{E}(X_{k,n})$ and $\sigma_{k,n}^2=\mathbb{V}ar(X_{k,n})$, $1\leq k \leq k(n)$, if these expectations or variances exist. We also suppose that 

\begin{equation*} 
k(n) \rightarrow +\infty \  as \ n\rightarrow +\infty. 
\end{equation*}

\Bin The central limit theorem problem consists in finding, whenever possible, the weak limit law (in type) of the by-row sums of the array $X$,  \textit{i.e.} the summands:

$$
S_n[X]=\sum_{k=1}^{k(n)} X_{k,n}, \ n\geq 1.
$$ 

\Bin Historically, the \textit{CLT} was discovered with the convergence of a Binomial law (which has the same law as a sum of \textit{iid} Bernoulli random variables) to the standard Gaussian law (due to Laplace, De Moivre, etc., around $1731$, see \cite{loeve} for a review). For a long period, the Gaussian limit was automatically meant in the \textit{CLT} problem. Many authors, among them L\'evy, Gnedenko, Kolmogorov, etc., characterized the class of possible limit laws under the \textit{Uniform Asymptotic Negligibility (\textit{UAN})} condition, exactly as the class of infinitely decomposable distributions. The longtime association of \textit{CLT}'s with Gaussian limits explains that some authors reserve the vocable \textit{CLT} for Gaussian limits and for other possible limits, they use different vocables. Here we use the vocable of $\textbf{$G$-\textit{CLT}}$ to cover all possible limit laws $G$ beyond the Gaussian law.\\

\Ni Here we suppose that the $X_{k,n}$'s are integrable with finite variances. For an array $X$, we define some important hypotheses used in the formulation of the \textit{CLT} problem.\\

\Ni \textbf{(1) The \textit{UAN} condition}: for any $\varepsilon>0$,

\begin{equation}
U(n,\varepsilon,X)= \sup_{1\leq k \leq k_n} \mathbb{P}(|X_{k,n} - a_{k,n}|\geq \varepsilon)\rightarrow 0.
\end{equation}

\Bin \textbf{(2) The Bounded Variance Hypothesis (\textit{BVH})}: there exists a constant $c>0$,

$$
\sup_{n\geq 1} MV(n,X) \leq c,
$$

\Bin where

$$
MV(n,X)=\mathbb{V}ar(S_n[X]), \ n\geq 1.
$$

\Bin \textbf{(3) The Variance Convergence Hypothesis (\textit{VCH})}:

$$
MV(n,X) \rightarrow c \in ]0,+\infty[.
$$

\Bin According to the state of the art in \textit{CLT}'s theory for centered, square integrable  and independent by-row arrays of random variables, the summands weakly converges to a probability law  associated to the \textit{cdf} $G$ and to the characteristic function (\textit{cha.f}) $\psi_{G}$ under the \textit{UAN} condition and the \textit{BVH} if and only if the sequence of distribution functions (\textit{df})

$$
K_n(x)=\sum_{k=1}^{k(n)} \int_{-\infty}^{x} y^2 dF_{k,n}(y), \ x\in \mathbb{R}, \ n\geq 1,
$$

\Bin pre-weakly converges to a \textit{df} $K$, denoted $K_n \rightsquigarrow_{pre} K$, that is for any continuity point $x$ of $K$ denoted as $[x \in C(K)]$, we have

$$
K_n(x) \rightarrow K(x),
$$

\Bin and the \textit{cha.f} $\psi_{G}(\circ)$ of $G$ is given by $\exp(\psi[K](\circ))$ with

\begin{equation*}
\forall u \in \mathbb{R}, \ \psi[K](u)=\int \frac{e^{iux}-1-iux}{x^2} \ dK(x).
\end{equation*}

\Bin If we have the \textit{VCH}, the convergence criterion is replaced by the weak convergence $K_n \rightsquigarrow K$. Moreover, the limit law $G$ is necessarily an infinitely decomposable law.\\

\Ni In the non centered case, with the same hypotheses above on the random variables of the array, the summands weakly converges to a probability law  associated to the \textit{cdf} $G^{\ast}$ and to the \textit{cha.f} $\psi_{G^{\ast}}$ under the \textit{UAN} condition and the \textit{BVH} if and only if 

$$
\sum_{k=1}^{k(n)} a_{k,n} \rightarrow a, \ a\in\mathbb{R}
$$

\Bin and the sequence of distribution functions (\textit{df})

$$
K_n^{\ast}(x)=\sum_{k=1}^{k(n)} \int_{-\infty}^{x} y^2 dF_{k,n}(y+a_{k,n}), \ x\in \mathbb{R}, \ n\geq 1,
$$

\Bin pre-weakly converges to a \textit{df} $K^{\ast}$ and the \textit{cha.f} $\psi_{G^{\ast}}(\circ)$ of $G^{\ast}$ is given by $\exp(\psi[K^{\ast}](\circ))$ with

\begin{equation*}
\forall u \in \mathbb{R}, \ \psi[K^{\ast}](u)=\int \frac{e^{iux}-1-iux}{x^2} \ dK^{\ast}(x).
\end{equation*}

\Bin If we have the \textit{VCH}, the convergence criterion is replaced by the weak convergence $K_n^{\ast} \rightsquigarrow K^{\ast}$. Moreover, the limit law $G^{\ast}$ is of the form $G^{\ast} = G+a$, with $G$ is necessarily a centered and infinitely decomposable law.\\

\Ni By specializing the limit law as a Gaussian law or a Poisson law, which clearly are infinitely decomposable laws, we have the following characterizations. \\

\Ni \textbf{C1}. Under the conditions

$$
\biggr(\forall n\geq 1, \ \forall 1\leq k\leq k(n), \ \ a_{k,n} = 0\biggr) \ \ and \ \ \sum_{k=1}^{k(n)}\sigma_{k,n}^2 = 1,
$$

\Bin the summands $S_n[X]$ of the array $X$ converges to standard Gaussian law and $\max_{1\leq k\leq k(n)}\sigma_{k,n}^2\rightarrow 0$ if and only if the following Lynderberg-Gaussian condition holds:

\begin{equation}
\forall \varepsilon>0, \ \ L_{n,G}(\varepsilon) = \sum_{k=1}^{k(n)} \int_{(|x|\geq \varepsilon)} x^2 \ dF_{k,n}(x)\rightarrow 0. \label{lyndGauss}
\end{equation}

\Bin \textbf{C2}. Under the conditions

$$
\max_{1\leq k\leq k(n)}\sigma_{k,n}^2\rightarrow 0 \ \ and \ \ \sum_{k=1}^{k(n)}\sigma_{k,n}^2 \rightarrow \lambda, \ \lambda>0,
$$

\Bin the summands $S_n[X]$ of the array $X$ converges to a  translated Poisson law $\mathcal{P}(a,\lambda)\equiv a+\mathcal{P}(\lambda)$, \ $a\in \mathbb{R}$,  if and only if 

$$
\sum_{k=1}^{k(n)} a_{k,n} \rightarrow a + \lambda
$$

\Bin and the following Lynderberg Poisson-type condition holds:

\begin{equation}
\forall \varepsilon>0, \ \ L_{n,P}(\varepsilon) = \sum_{k=1}^{k(n)} \int_{(|x-1|\geq \varepsilon)} x^2 \ dF_{k,n}(x+a_{k,n})\rightarrow 0. \label{lyndPoisson}
\end{equation}

\newpage
\section{Statements of the results} \label{sec-03}

\Bin As announced, we focus here on the non-stationary independent scheme. \\

\Bin First, we consider uniform conditions of the convergence of the probabilities $p_{k,n}$ (in the Bernoulli case) and $q_{k,n}$ (in the corrected geometric case) to zero to unveil refined versions of the extensions. Later, we will provide more general conditions. \\
 
\begin{theorem}\label{piibinomialIS} 

\Ni Let 

$$
X=\biggr\{ \{X_{k,n}, \ 1 \leq k \leq k_n=k(n)\}, \ n\geq 1\biggr\},
$$

\Bin be an array of by-row independent Bernoulli random variables, that is: \\

\Ni (1) $\forall n\geq 1$, $\forall 1\leq k \leq k(n)$, $X_{k,n}\sim \mathcal{B}(p_{k,n})$, with $0<p_{k,n}<1$ and:\\

\Ni (2) $\sup_{1\leq k \leq k(n)} p_{k,n} \rightarrow 0$;\\

\Ni (3) $\sum_{1\leq k \leq k(n)} p_{k,n}\rightarrow \lambda \in ]0, \ +\infty[$.\\

\Ni Then we have

$$
S_n[X] \rightsquigarrow \mathcal{P}(\lambda).
$$

\end{theorem}

\begin{theorem}\label{niibinomialIS}  

\Ni Let 

$$
X=\biggr\{ \{X_{k,n}, \ 1 \leq k \leq k_n=k(n)\}, \ n\geq 1\biggr\},
$$

\Bin be an array of by-row-independent corrected geometric random variables, that is: \\

\Ni (1) $\forall n\geq 1$, $\forall 1\leq k \leq k(n)$, $X_{k,n}\sim \mathcal{G}^{\ast}(p_{k,n})$, with $0<p_{k,n}=1-q_{k,n}<1$ and:\\

\Ni (2) $\sup_{1\leq k \leq k(n)} q_{k,n} \rightarrow 0$;\\

\Ni (3) $\sum_{1\leq k \leq k(n)} q_{k,n}\rightarrow \lambda \in ]0, \ +\infty[$.\\

\Ni Then we have

$$
S_n[X] \rightsquigarrow \mathcal{P}(\lambda).
$$
\end{theorem}

\Bin \textbf{A simple application}. Let us give a simple application to a classical example. We suppose that we observe occurrences of landing crashes at some airport (A) over a period $T>0$. We know that those crashes are usually of very low probabilities.  Over $n$ landings, we denote $X_n$ the number of crashes at times $k(n)$. Usually, we suppose that the data (of landing crashes) are observations of \textit{iid} Bernoulli $\mathcal{B}(p_n)$ and then, the approximation $X_n\approx Z\sim \mathcal{P}(\lambda)$, with $\lambda=np_n$ can be used. That formula was systematically used with limited performance of computers. However, with powerful computers, we no-longer need that approximation to compute the related p-values $\mathbb{P}(X_n>t)$ of the statistical tests since we know the explicit form of $S_n[X]$. In the software \textbf{R},  the code $1-pbiniom(t,p_n)$ gives the desired values.\\

\Bin Now suppose we can use the independence hypothesis only and not the stationary distribution. Hence the distribution of $X_n$ is the convolution product of  Bernoulli $\mathcal{B}(p_{k,n})$ distributions and its law is not simple. So, the simplest way to compute $\mathbb{P}(X_n>t)$ should be using the approximation  $\mathbb{P}(\mathcal{P}(\lambda)>t)$ with $\lambda=p_{1,n}+\cdots+p_{k(n),n}$. So it is better to use the non-stationary scheme since the stationary hypothesis is usually a \textit{working hypothesis}, not confirmed, and $p_n$ is computed as the average number of crashes.\\

\Bin These two theorems actually are still particular cases of two more general results.\\

\begin{theorem}\label{piibinomialISG} 

\Ni Let 

$$
X=\biggr\{ \{X_{k,n}, \ 1 \leq k \leq k_n=k(n)\}, \ n\geq 1\biggr\},
$$ 

\Bin be an array of by-row-independent Bernoulli random variables, that is: \\

\Ni (GP1) $\forall n\geq 1$, $\forall 1\leq k \leq k(n)$, $X_{k,n}\sim \mathcal{B}(p_{k,n})$, with $0<p_{k,n}<1$ and:\\

\Ni (GP2) $\sup_{1\leq k \leq k(n)} p_{k,n}(1-p_{k,n}) \rightarrow 0$;\\

\Ni (GP3) $\sum_{1\leq k \leq k(n)} p_{k,n}(1-p_{k,n}) \rightarrow \lambda \in ]0, \ +\infty[$ and $\sum_{1\leq k \leq k(n)} p_{k,n}  \rightarrow \lambda$;\\

\Ni (GP4) for $0<\varepsilon<1$, $n\geq 1$, $1\leq k\leq k(n)$, and for

$$
B(\varepsilon,k,n)=\sum_{j=0}^{1}  1_{\biggr(\biggr|j-q_{k,n}-1\biggr|\geq \varepsilon\biggr)} \biggr(j-p_{k,n}\biggr)^2 p_{k,n}^j q_{k,n}^{1-j},
$$

\Bin we have 

$$
B(\varepsilon,n)=\sum_{k=1}^{k(n)} B(\varepsilon,k,n) \rightarrow 0.
$$

\Bin Then we have

$$
S_n[X] \rightsquigarrow \mathcal{P}(\lambda).
$$

\end{theorem}

\begin{theorem}\label{niibinomialISG}  

\Ni Let 

$$
X=\biggr\{ \{X_{k,n}, \ 1 \leq k \leq k_n=k(n)\}, \ n\geq 1\biggr\},
$$ 

\Bin be an array of by-row-independent corrected geometric random variables, that is: \\

\Ni (GN1) $\forall n\geq 1$, $\forall 1\leq k \leq k(n)$, $X_{k,n}\sim \mathcal{G}^{\ast}(p_{k,n})$, with $0<p_{k,n}=1-q_{k,n}<1$ and:\\

\Ni (GN2) $\sup_{1\leq k \leq k(n)} (q_{k,n}/p_{k,n}^2) \rightarrow 0$;\\

\Ni (GN3) for $h \in \{1,2\}$, \ \ $\sum_{1\leq k \leq k(n)} (q_{k,n}/p_{k,n}^h)\rightarrow \lambda \in ]0, \ +\infty[$.\\

\Ni (GN4) for $0<\varepsilon<1$, $n\geq 1$, $0\leq k\leq k(n)$, and for

$$
B(\varepsilon,k,n)=\sum_{j=0}^{+\infty}  1_{\biggr(\biggr|j-\frac{q_{k,n}}{p_{k,n}}-1\biggr|\geq \varepsilon\biggr)} \biggr(j-\frac{q_{k,n}}{p_{k,n}}\biggr)^2 p_{k,n} q_{k,n}^j,
$$

\Bin we have

$$
B(\varepsilon,n)=\sum_{k=1}^{k(n)} B(\varepsilon,k,n) \rightarrow 0.
$$

\Ni Then we have

$$
S_n[X] \rightsquigarrow \mathcal{P}(\lambda).
$$
\end{theorem}

\Bin \textbf{Proofs of Theorems \ref{piibinomialIS} and \ref{niibinomialIS}}. Let us begin by:\\

\Ni \textit{(A1)- Proof of Theorem \ref{piibinomialIS}}. Throughout this proof, the notation $\ell_{k,n}=\overline{o}_n(1)$, for $k$ ranging over some set $I_n$ means that the sequence $\ell_{k,n}$ goes to zero as $n\rightarrow +\infty$ uniformly in $k \in I_n$. So Assumption (2) means that

$$
p_{k,n}=\overline{o}_n(1) \ \ and \ \ q_{k,n}=1-p_{k,n}=1+\overline{o}_n(1).
$$

\Bin We have to check the \textit{\textbf{UAN}} condition. By using Chebychev's inequality, we have, for any $\varepsilon>0$,
\label{GPU}
\begin{eqnarray*}
U(n,\varepsilon,X)&=&\sup_{1\leq k \leq k_n} \mathbb{P}(|X_{k,n} - a_{k,n}|\geq \varepsilon)\\
&\leq& \varepsilon^{-2} \sup_{1\leq k \leq k_n} \mathbb{V}ar(X_{k,n})\\
&=& \varepsilon^{-2} \ p_{k,n} \ q_{k,n} \\
&=&\varepsilon^{-2} \ \overline{o}_n(1) (1+\overline{o}_n(1))\rightarrow 0.
\end{eqnarray*}

\Bin The \textit{\textbf{VCH}} also holds since

\begin{eqnarray*}
MV(n,X)&=&\sum_{1\leq k \leq k(n)} \mathbb{V}ar(X_{k,n})\\
&=&\sum_{1\leq k \leq k(n)} p_{k,n} \ q_{k,n}\\
&=& (1+\overline{o}_n(1)) \sum_{1\leq k \leq k(n)} p_{k,n}\rightarrow \lambda.
\end{eqnarray*}

\Bin Besides

$$
\sum_{1\leq k \leq k(n)} \mathbb{E}(X_{k,n})=\sum_{1\leq k \leq k(n)} p_{k,n}\rightarrow \lambda.
$$

\Bin So, we are in the position of applying the conditions of weak convergence to a Poisson law by checking the Poisson Lynderbeg condition \eqref{lyndPoisson}. We have
for any $\varepsilon>0$

$$
L_{n,P}(\varepsilon) =: \sum_{k=1}^{k(n)} L_{n,k,P}(\varepsilon),
$$

\Bin with 

\begin{eqnarray*}
L_{n,k,P}(\varepsilon)&=&\int_{(|x-1|\geq \varepsilon)} x^2 \ dF_{k,n}(x+p_{k,n})\\
&=& \int_{(|X_{k,n}-p_{k,n}-1|\geq \varepsilon)} |X_{k,n}-p_{k,n}|^2 \ d\mathbb{P}_{k,n}\\
&=& p_{k,n}\biggr(1_{(|X_{k,n}-p_{k,n}-1|\geq \varepsilon)} |X_{k,n}-p_{k,n}|^2\biggr)_{(X_{k,n}=1)}\\
&+& (1-p_{k,n})\biggr(1_{(|X_{k,n}-p_{k,n}-1|\geq \varepsilon)} |X_{k,n}-p_{k,n}|^2\biggr)_{(X_{k,n}=0)}\\
&=&p_{k,n} 1_{(|p_{k,n}|\geq \varepsilon)} \ (1-p_{k,n})^2 + (1-p_{k,n})1_{(|p_{k,n}+1|\geq \varepsilon)} \ p_{k,n}^2\\
&=&p_{k,n} 1_{(|\overline{o}_n(1)|\geq \varepsilon)} (1+\overline{o}_n(1))^2 + \overline{o}_n(1)(1+\overline{o}_n(1)) 1_{(|\overline{o}_n(1)+1|\geq \varepsilon)} \ p_{k,n}.
\end{eqnarray*}

\Bin We only need to get \eqref{lyndPoisson} for $0<\varepsilon<\varepsilon_0$, for a fixed $\varepsilon_0>0$. Let us fix $\varepsilon_0=1/2$. So, for $n$ large enough, 

$$
1_{(|\overline{o}_n(1)|\geq \varepsilon)}=0 \ \ and \ \ 1_{(|\overline{o}_n(1)+1|\geq \varepsilon)}=1
$$

\Bin and hence

\begin{eqnarray*}
\sum_{k=1}^{k(n)} L_{n,k,P}(\varepsilon) &=&  \overline{o}_n(1) (1+\overline{o}_n(1)) \sum_{k=1}^{k(n)} p_{k,n}\\
&=& \overline{o}_n(1) (1+\overline{o}_n(1)) (\lambda + o(1))\rightarrow 0.
\end{eqnarray*}

\Bin The proof is complete.  $\blacksquare$\\

\newpage
\Ni \textit{(A2)- Proof of Theorem \ref{niibinomialIS}}. Assumption (2) of the theorem means that

$$
q_{k,n}=\overline{o}_n(1), \  p_{k,n}=1+\overline{o}_n(1) \ \ and \ \ (1/p_{k,n})^i=1+\overline{o}_n(1), \ i=1,2,3.
$$

\Bin We have to check the \textit{\textbf{UAN}} condition. By using Chebychev's inequality, we have, for any $\varepsilon>0$,

\begin{eqnarray*}
U(n,\varepsilon,X)&=&\sup_{1\leq k \leq k_n} \mathbb{P}(|X_{k,n} - a_{k,n}|\geq \varepsilon)\\
&\leq& \varepsilon^{-2} \sup_{1\leq k \leq k_n} \mathbb{V}ar(X_{k,n})\\
&=& \varepsilon^{-2} \sup_{1\leq k \leq k_n} \frac{q_{k,n}}{p_{k,n}^2}\\
&=&\varepsilon^{-2} \overline{o}_n(1) (1+\overline{o}_n(1))\rightarrow 0.
\end{eqnarray*}

\Bin The \textit{\textbf{VCH}} also holds since

\begin{eqnarray*}
MV(n,X)&=&\sum_{1\leq k \leq k(n)} \mathbb{V}ar(X_{k,n})\\
&=&\sum_{1\leq k \leq k(n)} \frac{q_{k,n}}{p_{k,n}^2}\\
&=& (1+\overline{o}_n(1)) \sum_{1\leq k \leq k(n)} q_{k,n}\rightarrow \lambda.
\end{eqnarray*}

\Bin Besides

\begin{eqnarray*}
\sum_{1\leq k \leq k(n)} \mathbb{E}(X_{k,n})=\sum_{1\leq k \leq k(n)} \frac{q_{k,n}}{p_{k,n}}
&=&(1+\overline{o}_n(1)) \sum_{1\leq k \leq k(n)} q_{k,n} \rightarrow \lambda.
\end{eqnarray*}

\Bin Here again, we are in the position of applying the conditions of weak convergence to a Poisson law by checking the Poisson Lynderbeg condition \eqref{lyndPoisson}. We have
for any $0<\varepsilon<1/2$

$$
L_{n,P}(\varepsilon) =: \sum_{k=1}^{k(n)} L_{n,k,P}(\varepsilon),
$$

\Bin with 

\begin{eqnarray*}
L_{n,k,P}(\varepsilon) &=&\int_{(|x-1|\geq \varepsilon)} x^2 \ dF_{k,n}(x+a_{k,n})\\
&=& \int_{\biggr(\left|X_{k,n}-\frac{q_{k,n}}{p_{k,n}}-1\right|\geq \varepsilon\biggr)} \left|X_{k,n}-\frac{q_{k,n}}{p_{k,n}}\right|^2 \ d\mathbb{P}_{X_{k,n}}\\
&=& \sum_{j=0}^{+\infty} p_{k,n} q_{k,n}^j \biggr(1_{\left(\left|X_{k,n}-\frac{q_{k,n}}{p_{k,n}}-1\right|\geq \varepsilon \right)} \left|X_{k,n}-\frac{q_{k,n}}{p_{k,n}}\right|^2 \biggr)_{(X_{k,n}=j)}\\
&=& p_{k,n}  1_{\left(\left|\frac{q_{k,n}}{p_{k,n}}+1\right|\geq \varepsilon \right)} \left(\frac{q_{k,n}}{p_{k,n}}\right)^2  \ \ (for \ j=0)\\
&+& p_{k,n} q_{k,n}  1_{\left(\left|\frac{q_{k,n}}{p_{k,n}}\right|\geq \varepsilon \right)} \left(1-\frac{q_{k,n}}{p_{k,n}}\right)^2  \ \ (for \ j=1)\\
&+& \sum_{j=2}^{+\infty} p_{k,n} q_{k,n}^j  1_{\left(\left|j-1-\frac{q_{k,n}}{p_{k,n}}\right|\geq \varepsilon \right)} \left(j-\frac{q_{k,n}}{p_{k,n}}\right)^2  \ \ (for \ j\geq 2).\\
\end{eqnarray*}

\Bin By the same remarks used in the precedent proof, we have for $n$ large enough, 

$$
1_{\left(\left|\frac{q_{k,n}}{p_{k,n}}+1\right|\geq \varepsilon \right)}= 1_{\left(\left|\overline{o}_n(1)(1+\overline{o}_n(1))+1\right|\geq \varepsilon \right)}=1
$$

\Bin and next

$$
p_{k,n}  1_{\left(\left|\frac{q_{k,n}}{p_{k,n}}+1\right|\geq \varepsilon \right)} \left|\frac{q_{k,n}}{p_{k,n}}\right|^2=(1+\overline{o}_n(1)) q_{k,n}^2 = \overline{o}_n(1)(1+\overline{o}_n(1)) \ q_{k,n}. \ \ (L1) 
$$

\Bin Also, we have for $n$ large enough,

$$
1_{\left(\left|\frac{q_{k,n}}{p_{k,n}}\right|\geq \varepsilon \right)}=1_{\left(\left|\overline{o}_n(1)(1+\overline{o}_n(1))\right|\geq \varepsilon \right)}=0
$$

\Bin and next
$$
p_{k,n} q_{k,n}  1_{\left( \left|\frac{q_{k,n}}{p_{k,n}}\right|\geq \varepsilon\right)} \left(1-\frac{q_{k,n}}{p_{k,n}}\right)^2=0. \ \ (L2)
$$

\Bin Now, for $n$ large enough and for any $j\geq 2$,

$$
1_{(|j-1+\overline{o}_n(1)(1+\overline{o}_n(1))|\geq \varepsilon)}=1
$$

\Bin and thus,

\begin{eqnarray*}
A_{j,k,n}&:=&p_{k,n} q_{k,n}^j  1_{\biggr(\biggr|j-1-\frac{q_{k,n}}{p_{k,n}}\biggr|\geq \varepsilon\biggr)} \left(j-\frac{q_{k,n}}{p_{k,n}}\right)^2 \notag\\
&=& \biggr((1+\overline{o}_n(1))q_{k,n}\biggr) \ q_{k,n}^{j-1}  1_{\biggr(\biggr|j-1+\overline{o}_n(1)(1+\overline{o}_n(1))\biggr|\geq \varepsilon\biggr)} \biggr(j+\overline{o}_n(1)(1+\overline{o}_n(1))\biggr)^2 \notag\\
&=& \biggr((1+\overline{o}_n(1))q_{k,n}\biggr)\ q_{k,n}^{j-1} (j+\overline{o}_n(1)(1+\overline{o}_n(1)))^2 \notag\\
&\leq & \biggr((1+\overline{o}_n(1))q_{k,n}\biggr)\ q_{k,n}^{j-1}   (j+1)^2\notag\\
&\leq & 2\biggr((1+\overline{o}_n(1))q_{k,n}\biggr)\ q_{k,n}^{j-1}   (j^2+1),
\end{eqnarray*}

\Bin where we apply the $C_2$-inequality in the last line. So we have

\begin{eqnarray*}
\sum_{j\geq 2} A_{j,k,n}&\leq& 2(1+\overline{o}_n(1)) \ q_{k,n}\ B(n,k), \ \ (L3)
\end{eqnarray*}

\Bin with

$$
B(n,k)=\sum_{j\geq 2} q_{k,n}^{j-1} + \sum_{j\geq 2} j^2q_{k,n}^{j-1}=:B(n,k,1)+B(n,k,2).
$$

\Bin We have

$$
B(n,k,1)=\biggr(\sum_{j\geq 0} q_{k,n}^{j}\biggr)-1=\frac{q_{k,n}}{p_{k,n}} = \overline{o}_n(1)(1+\overline{o}_n(1)). \ \ (L4a)
$$

\Bin Next

\begin{eqnarray*}
B(n,k,2)&=&\sum_{j\geq 2} j q_{k,n}^{j-1} + \sum_{j\geq 2} j(j-1) q_{k,n}^{j-1}\\
&=& \biggr(\sum_{j=1}^{+\infty} j q_{k,n}^{j-1} -1\biggr)+\biggr(q_{k,n}\sum_{j=2}^{+\infty} j(j-1) q_{k,n}^{j-2}\biggr)\\
&=&\biggr(\biggr\{ \sum_{j=0}^{+\infty} q_{k,n}^{j}\biggr\}^{\prime} -1\biggr)+q_{k,n}\biggr\{\sum_{j=0}^{+\infty} q_{k,n}^{j}\biggr\}^{\prime\prime}\\
&=&\biggr(\frac{1}{p_{k,n}^2} -1\biggr)+\biggr(q_{k,n}\frac{2}{p_{k,n}^3}\biggr)\\
&=&\frac{1-p_{k,n}^2 }{p_{k,n}^2} + \frac{2q_{k,n}}{p _{k,n}^3} \\
&=&\frac{p_{k,n}(1-p_{k,n}^2) + 2q_{k,n}}{p_{k,n}^3} \\
&=&\frac{p_{k,n}q_{k,n}(1+p_{k,n}) + 2q_{k,n}}{p_{k,n}^3} \\
&=&\frac{q_{k,n}\left(p_{k,n}(1+p_{k,n}) + 2\right)}{p_{k,n}^3} \\
&\leq& \frac{4q_{k,n}}{p_{k,n}^3} \\
&=& 4 \overline{o}_n(1)(1+\overline{o}_n(1)). \ \ (L4b)
\end{eqnarray*}

\Bin Hence 

$$
B(n,k) \leq C \overline{o}_n(1)(1+\overline{o}_n(1)), \ \  (L4c)
$$

\Bin for some $C>0$ by (L4a) and (L4b). \\

\Ni Finally, by putting together (L1), (L2), (L3) and (L4c), we get

$$
L_{n,P}(\varepsilon)\leq \overline{o}_n(1)(1+\overline{o}_n(1))(\lambda + o(1)) \biggr\{1+2 C(1+\overline{o}_n(1))\biggr\}\rightarrow 0.
$$

\Bin This completes the proof.  $\blacksquare$\\

\Ni \textbf{Proofs of Theorems \ref{piibinomialISG} and \ref{niibinomialISG}}. \\

\Ni In the proofs of Theorems \ref{piibinomialIS} and \ref{niibinomialIS}, the \textit{\textbf{UAN}}, the \textit{\textbf{CVH}}, the Poisson Lynderberg condition and the convergence of $\mathbb{E}S_n[X]$ are the general conditions for $S_n[X]\rightsquigarrow \mathcal{P}(\lambda)$. In Theorems \ref{piibinomialISG} and \ref{niibinomialISG}, we used uniformity convergence to simplify these conditions. $\square$
\newpage

\section{Concluding remarks} \label{sec-04} 

\Ni The extensions we provide are the first general results. The central limit theorem frame seems to be the appropriate way to get more the general extensions. Theorems \ref{piibinomialIS} and \ref{niibinomialIS} can be done by direct methods. However, general forms in Theorems \ref{piibinomialISG} and \ref{niibinomialISG} could hardly be obtained in direct methods. They are products of the \textit{CLT} frame.

\newpage 
\Ni \textbf{Appendix : Proof of Proposition \ref{niidBinomial}}.  \label{niidBinomial-proof}\\

\Ni Let us use the convergence of characteristic functions. Let $X_{n}$ be a sequence of $\mathcal{NB}(n,p_n)$-random variables and $X$ be a $\mathcal{P}(\lambda)$ random variable. We have 

$$
S_n=X_n-n=Z_1+\cdots+Z_n,
$$

\Bin where $Z, \ Z_1, \cdots, \ Z_n$ are independent random variables such that each $Z_i+1$ follows a geometric law of parameter $p_n$, that is

$$
\Phi_{Z}(t)=\frac{p_n}{1-q_ne^{it}}, \ \ t\in \mathbb{R}.
$$

\Bin So, for $t\in \mathbb{R}$ fixed,

\begin{eqnarray*}
\Phi _{S_{n}}(t)=\exp\biggr(n \left(\log p_n - \log\left(1-q_ne^{it}\right)\right)\biggr).
\end{eqnarray*}

\bigskip \noindent We have, as $n\rightarrow +\infty$,

$$
n \log p_n= n \log (1-q_n)=-nq_n +o(nq_n)
$$

\noindent and

$$
-n \log\left(1-q_ne^{it}\right) = nq_n e^{it} + o(nq_n).
$$

\Bin Hence we get for any $t\in \mathbb{R}$,

\begin{eqnarray*}
\Phi_{S_{n}}(t)=\exp\biggr(nq_n (e^{it}-1) + o(nq_n)\biggr) \rightarrow e^{\lambda(e^{it}-1)}=\Phi_{\mathcal{P}(\lambda)}(t). \ \square
\end{eqnarray*}

\end{document}